\documentclass[12pt]{amsart}
\usepackage[all]{xy}
\usepackage{amssymb}
\usepackage{amsmath}
\usepackage{MnSymbol}
\usepackage{setspace}
\usepackage{graphicx}
\usepackage{graphics,graphicx}
\usepackage{pgf,color}
\usepackage{epstopdf}
\usepackage{epsfig,textpos,mathrsfs}
\usepackage{amsthm}
\newtheorem{definition}{Definition}[section]
\usepackage[utf8]{inputenc}
\usepackage[english]{babel}
\usepackage{amsthm}
\newtheorem{remark}{Remark}[section]
\newtheorem{theorem}{Theorem}[section]
\newtheorem{lemma}{Lemma}[section]
\newtheorem{proposition}{Proposition}[section]
\newtheorem{corollary}{Corollary}[section]
\newtheorem{example}{Example}[section]

\theoremstyle{remark}

\newcommand{\bt}{\begin{theorem}}
\newcommand{\et}{\end{theorem}}
\newcommand{\bl}{\begin{lemma}}
\newcommand{\el}{\end{lemma}}
\newcommand{\bexc}{\begin{exercise}}
\newcommand{\eexc}{\end{exercise}}
\newcommand{\bp}{\begin{proposition}}
\newcommand{\ep}{\end{proposition}}
\newcommand{\bex}{\begin{example}}
\newcommand{\eex}{\end{example}}
\newcommand{\bc}{\begin{corollary}}
\newcommand{\ec}{\end{corollary}}
\newcommand{\bo}{\begin{proof}}
\newcommand{\eo}{\end{proof}}
\newcommand{\bd}{\begin{definition}}
\newcommand{\ed}{\end{definition}}
\newcommand{\br}{\begin{remark}}
\newcommand{\er}{\end{remark}}
\newcommand{\be}{\begin{enumerate}}
\newcommand{\ee}{\end{enumerate}}

\newcommand{\T}{\mathbb{T}}

\newcommand{\A}{{\mathcal A}}

\newcommand{\cO}{{\mathcal O}}

\newcommand{\Z}{{\mathbb Z}}

\newcommand{\N}{{\mathbb N}}

\newcommand{\C}{{\mathbb C}}

\begin{document}

\title{Rendezvous with Sensitivity}
\author{Anima Nagar}
\address{Department of Mathematics, Indian Institute of Technology Delhi,\\
Hauz Khas, New Delhi 110016, INDIA}
\email{anima@maths.iitd.ac.in}

\maketitle

\renewcommand{\thefootnote}{}



\footnote{This article is dedicated to Joseph Auslander and James Yorke.}

\footnote{ This article is based on a talk given at the International Workshop and Conference on Topology \& Applications, Cochin, INDIA, Dec. 2018.}

\renewcommand{\thefootnote}{\arabic{footnote}}
\setcounter{footnote}{0}

\begin{abstract} Let $(X,d)$ be a compact metric space and $f:X \to X$ be a self-map.
The compact dynamical system $(X,f)$ is called sensitive or sensitivity depends on initial conditions, if there is a positive constant $\delta$ such that in each non-empty open subset there are  distinct points
whose iterates will be $\delta-$apart  at  same instance. This dynamical property, though being a very weak one, brings in the essence of unpredictability in the system. In this article, we survey various sensitivities and some properties implied by and implying such sensitivities. \end{abstract}

\begin{section}{Introduction}

\emph{Sensitive Dependence on Initial Conditions} is widely understood as
being the central idea of \emph{``chaos''} and was popularized by the meteorologist
\textbf{Edward Lorenz} through the so-called \emph{``butterfly effect''}(1963). The word \emph{``sensitive dependence on initial conditions''} was first used by \textbf{D. Ruelle} in 1977.

The foundational contributions to the general theory of \emph{``chaos''} can be credited to \textbf{Poincare, Birkhoff, Landau, Ruelle, Takens, Lorenz, Rossler, Bowen, May,  Liapunov, Sinai, Sharkovsky and many others who have not been named here}.

The  systematization of the mathematical theory and many initial mathematical definitions of \emph{``Chaos''} can be credited to  \textbf{James Yorke} \cite{ly}.  Another important concept has been the virtue of systems with every orbit dense in the space. This notion of \emph{``minimality''}  has had its widest exploration by \textbf{Joseph Auslander} \cite{aus}.

We look into the concept of \emph{``sensitivity''} which largely combines the view of \emph{``chaotic sets''} due to \textbf{Yorke} combined with  the notion of \emph{``minimality''} due to \textbf{Auslander}.

We note here that at the same time as when \textbf{Auslander and Yorke} were into publishing their work on sensitivity, there was defined a parallel definition of sensitivity by \textbf{John Guckenheimer} in a very restricted sense. His approach was mainly restricted to  special kind of interval maps \cite{gu}. \textbf{We follow the general definition due to Auslander and Yorke.}

\end{section}

\medskip

\begin{section}{ Compact Dynamical Systems}

\bd
A pair $(X,f)$ is called a \textbf{compact dynamical system}, where
 $(X,d)$ is a compact metric space and $f: X \to X$ is a continuous surjection.
\ed

\medskip

Here   $\cO(x) = \{f^n(x) : n \in \N\}$ is called the \textbf{orbit} of the point $x$.

\bd For a  system $(X,f)$, a point $x_0 \in X$ is called a \textbf{fixed point} if $f(x_0) = x_0$. And $y_0 \in X$ is called a \textbf{periodic point} if there exists $n \in \N$ such that $f^n(y_0) = y_0$. The smallest such $n$ is called the \textbf{period} of $y_0$. \ed

\bd The \textbf{$\omega$-limit set} of a point $x \in X$ under $f$, denoted as
$\omega(x)$, is the set of all limit points of
$\cO(x)$. \ed

\bd    A point $x \in X$ is said to be
\textbf{non-wandering} if for every neighbourhood $U$ of $x$ there is a $n
\in \mathbb{N}$ such that $f^n(U) \cap U \neq \emptyset$. The set of all
non-wandering points of $f$ is denoted as $\Omega(f)$. \ed

\bd For a system $(X,f)$,  a point $x \in X$ is called \textbf{recurrent}   when $x \in \omega (x)$, i. e. for every open set $U$ containing $x$ there exist $j \in \N$ such that $f^j(x) \in U$. Equivalently there is a sequence $n_k \nearrow \infty$ such that $f^{n_k}(x) \to x$.\ed

\medskip

\bd A subset $E\subset X$ is called \textbf{$f$-invariant} if $f(E)\subseteq E$.\ed

Note that $\omega(x)$ is a closed, non-empty, $f$-invariant set with $f(\omega(x)) = \omega(x)$.

\medskip

\bd A \textbf{minimal set} is a non-empty closed invariant set, which
contains no proper non-empty closed invariant subsets.

Equivalently,
a nonempty set $M \subset X$ is minimal if, for each $x \in M$, the orbit closure $\overline{\cO(x)} = M$. \ed

Fixed points and periodic orbits are minimal
sets. Zorn's lemma shows that  in a compact dynamical system  $(X, f)$ any nonempty closed $f$-invariant subsets of $X$ contain minimal subsets.

\bd A point $x \in X$ is called an \textbf{almost periodic point} if its orbit closure $\overline{\cO(x)}$ is a minimal set. \ed

Note that an almost periodic point is always recurrent, and periodic points are always almost periodic.

\bd A compact system $(X,f)$ is
called  \textbf{topologically transitive} if for every pair of nonempty
open sets $U,V$ in $X$, there is a $n \in \mathbb{N}$
such that $f^n(U) \cap V \neq \emptyset$ [Equivalently,  $U \cap f^{-n}(V)  \neq \emptyset$]. \ed

\bp For the compact system $(X, f)$, the following are
equivalent:

( i ) $(X, f)$ is topologically transitive.

(ii) There is an $x \in X$ with a dense orbit.

(iii) There is an $x \in X$ with $\omega(x) = X$.\ep

A point with a dense orbit is called a \textbf{transitive point}.

Note that a compact system is minimal when each of its points is a transitive point.

\bd Let $(X,f)$ and $(Y,g)$ be compact systems.  A continuous surjective map
$\phi \colon X \to Y$ is called a \textbf{equivariant map} if $\phi(fx) = g\phi(x)$ for all
$x\in X$.

In addition if $\phi$ is a homeomorphism, we call $\phi$ a \textbf{conjugacy} and
say that $(X,f)$ and $(Y,g)$ are \textbf{conjugate systems}.\ed

Note that

\be
\item If $x \in X$ is a periodic point, then $\phi(x)$ is a periodic point in $Y$.
\item If $x \in X$ is a transitive point, then $\phi(x)$ is a transitive point in $Y$.
\item   If $y \in Y$ is a periodic point and $\phi$ is a conjugacy then  $\phi^{-1}(y)$ is a periodic point in $X$.
\item   If $y \in Y$ is a transitive point and $\phi$ is a conjugacy then  $\phi^{-1}(y)$ is a transitive point in $X$.
\ee

Note that for a periodic point $y \in Y$, and an equivariant map $\phi \colon X \to Y$, $\phi^{-1}(y)$ need not contain any periodic point in $X$. Same holds for transitive points. A simple example can be a union of minimal systems collapsing to a point, however we can say something concrete for almost periodic points.

\bp Let $\phi \colon X \to Y$ be a continuous surjection where $(X,f)$ and $(Y,g)$ are compact systems. If $y_0 \in Y$ is almost periodic then $\phi^{-1}(y_0) \subset X$ contains an almost periodic point in $X$. \ep

\bo Let $y_0 \in Y$ is almost periodic. Then $M_0 = \overline{\cO(y_0)} \subset Y$ is minimal. $\phi^{-1}(M_0)$ is a closed nonempty $f$-invariant subset of $X$ and so contains a minimal subset $N_0$. Now $\phi(N_0) = M_0$ and so there exists $x_0 \in N_0$ such that $\phi(x_0) = y_0$. Thus $\phi^{-1}(y_0) \subset X$ contains an almost periodic point. \eo

\end{section}

\begin{section}{Equicontinuity and Sensitivity}

\bd Let $(X, f)$ be a compact dynamical system. A point $x \in X$ is said to
be \textbf{equicontinuous(stable)} if, for any $\epsilon > 0$, there is a $\delta > 0$ such that $d(x, x') < \delta$ implies $d(T^n(x), T^n(x')) < \epsilon$, for $n = 0, 1, 2, \ldots.$

If all $x \in X$ are equicontinuous, then the system \textbf{$(X,f)$ is called equicontinuous}.\ed

\bd Thus $x \in X$ is \textbf{sensitive(not stable / unstable)} if and only if there is a $\delta > 0$ such that there is a sequence $\{x_j\}$ in $X$ and $n \in \N$ such that $x_j \to x$ and $d(f^n(x), f^n(x_j)) > \delta$.

Equivalently, $(X,f)$ is \textbf{sensitive} at $x \in X$ if there is a $\delta > 0$ such that for any $\epsilon > 0$ there exists  $y \in X$ with $d(x,y)< \epsilon$  and $n \in \N$ with  $d(f^n(x), f^n(y)) > \delta$.

 And if all $x \in X$ are sensitive, then the system \textbf{$(X,f)$ is called sensitive or sensitively dependent on initial conditions.}

\ed

\br Note that the uniform $\delta$ will follow due to compactness. \er

We note that if $x \in X$ is sensitive then there will be infinitely many such $n \in \N$ for which the above holds.

\bp \cite{ay} Let $(X, f)$ be a compact system.

( i ) If $x \in X$ is $\epsilon-$sensitive, then every point of $\overline{\cO(x)}$ is $\epsilon/2$ sensitive.

(ii) If $(X, f)$ is minimal and contains a sensitive point, then there
is an $\epsilon > 0$ such that every point is $\epsilon-$sensitive. \ep

We note that this follows easily by triangular inequality and the fact that  points in the orbit of  sensitive points are sensitive.

\br The above gives a nice dichotomy result - all minimal systems are either equicontinuous or sensitive \cite{ay}. This is generally called the \textbf{Auslander-Yorke dichotomy.} \er

We look into some examples:

\bex  \label{tm} Let $X= [0,1]$, the unit interval. Define $f: X \to X$ as
$$f(x) = \left\{
                                                                          \begin{array}{ll}
                                                                            2x, & \hbox{$0 \leq x \leq 1/2$;} \\
                                                                            2(1-x), & \hbox{$1/2 \leq x \leq 1$.}
                                                                          \end{array}
                                                                        \right.$$
 Then for the system
$(X,f)$, if $|x - y| < \frac{1}{2^{n+1}}$, then $|f^n (x) - f^n(y)| = 2^n|x -y| < \frac{1}{2}$ giving a sensitive system.
\eex

\bex \label{sm}   Let $\A$ be a finite set and define $X=
\A ^\N$.  Equip $X$ with product topology.  Then $X$ is a compact
metrizable space.  One of the compatible metrics on $X$ is:

$$d(x,y) = \inf \{ {1\over k+1} \mid x(n) = y(n) ~~{\rm for }~~ |n
|<k \}.$$  Define $\sigma: X \to X$ by $\sigma(x)(n) = x(n+1)$.  Then $\sigma$
is a continuous self map of $X$ -   giving a compact system $(X, \sigma)$.

For any configuration $x = x_0x_1x_2 \ldots \in X$, and $\epsilon > 0$ there exists a configuration $y = y_0y_1y_2 \ldots \in X$ such that $d(x,y) < \epsilon$ but there is an $n \in \N$ for which $d(\sigma^n(x), \sigma^n(y)) > 1/2$ - giving a sensitive system. \eex

\bex \label{ms} We recall Example \ref{sm}. Let $\Lambda= \{0,1\}$  and define $X=
\Lambda ^\Z$. We consider the shift map $\sigma: X \to X$.

We obtain a minimal subset of $X$, by constructing an almost periodic point $p \in X$ since then ${\overline{\cO(p)}}$ will be minimal. We take the classical construction due to Marston Morse and   Axel Thue, giving the Morse-Thue sequence.

This construction is done using substitution: $0 \to 01, \ 1 \to 10$. Hence,
$$0 \to 01 \to 0110 \to 01101001 \to 0110100110010110 \to \cdots$$
This will finally converge to some $x \in \{0,1\}^\N$. This construction indicates that every finite word in $x$ occurs syndetically often. Extend $x$ to $p \in X$ by defining $p(n) = \left\{
                                                                  \begin{array}{ll}
                                                                    x(n), & \hbox{$n \geq 1$;} \\
                                                                    x(-n-1), & \hbox{$n < 0$.}
                                                                  \end{array}
                                                                \right.$

Every word in $p$ occurs syndetically and $p$ is symmetric at the mid point, and so $p$ is almost periodic. Thus $({\overline{\cO(p)}}, \sigma)$ is a minimal dynamical system.

It is easy to note that the system $({\overline{\cO(p)}}, \sigma)$ is sensitive.
\eex

\bex  \label{ir} Let  $X= \T =\{z \in \C: |z|=1\}$ be the unit circle.  And let $g: X \to X$ be defined as $g(\theta) = \theta + 2 \pi \alpha (\mod 2 \pi)$, where $\alpha$ is an irrational number. Then $g$ is an isometry and  $(X,g)$ is a minimal  equicontinuous system.\eex

\end{section}

\begin{section}{Sensitivity and Chaos}

As discussed earlier, sensitivity plays an important role in ``Chaos''. The (possibly)second mathematical definition of chaos is given as,

\bd  \cite{ay}[\textbf{Auslander and Yorke Chaos}] The compact system $(X, f)$ is said to be \textbf{ chaotic} if it is \emph{``topological transitive''} and  \emph{``pointwise sensitive''}.\ed

Since ``pointwise sensitive'' is equivalent to ``sensitivity'' and ``topological transitivity'' is equivalent to the ``existence of dense orbit'', for a compact system \emph{Auslander-Yorke Chaos} can be defined as a sensitive system with a dense orbit.

\br We note that this definition of Chaos makes sense even for non-compact systems. \er

Historically, this definition of Chaos did not catch any attention. It was a modification of this definition given a few years later that caught the fancy of  the entire mathematical fraternity.

\bd \cite{de} [\textbf{Devaney Chaos}] A compact system $(X,f)$ is called \textbf{Devaney Chaotic} if:

1. $(X,f)$ is topologically transitive.

2. $(X,f)$ is sensitive.

3. $(X,f)$ admits a dense set of periodic points. \ed

Thus by adding the condition of dense set of periodic points, Auslander-Yorke Chaos became Devaney Chaos and a well accepted definition too. But all this came with a price. Since a few years later a redundancy was reported in this definition. It was deduced
that transitivity along with the denseness of periodic orbits
actually implies sensitivity, giving a redundancy  proved in \cite{bb, sil, gw} - where in
\cite{bb, sil} the space is not even assumed to be compact.

Infact, there are various conditions under which transitivity implies sensitivity. We refer to \cite{nk} for a survey on this.

\medskip

Despite this redundancy, sensitivity continues to be an integral part of chaos. Proving a system to be sensitive is much easier than proving it to be Devaney Chaotic.

\medskip

It becomes meaningful at this juncture to talk about the  first mathematical definition of Chaos given  by Li and Yorke \cite{ly}, which did generate a lot of excitement among the mathematical fraternity.

\bd For the system $(X,f)$, the pair  $(x, y) \in X \times X$ is called \textbf{proximal} if $\liminf \limits_n d(f^n(x), f^n(y)) = 0$, \textbf{asymptotic} if $\lim \limits_n d(f^n(x), f^n(y)) = 0$, \textbf{not $\delta-$asymptotic}  if $\limsup \limits_n d(f^n(x), f^n(y))  > \delta$, \textbf{Li-Yorke} if $(x,y)$ is proximal but not asymptotic and \textbf{$\delta-$Li-Yorke} if $(x,y)$ is proximal but not $\delta-$asymptotic.\ed

\bd \cite{ly}  A set $S \subseteq  X$ containing at least two points is called a \textbf{scrambled set} if for any $x, y \in S$ with $x \neq y$, the pair $(x,y)$ is Li-Yorke.  \ed

This means that  for any two distinct $x,y \in S$, the orbits of $x $ and $y$ get arbitrarily close to each other but infinitely many times they are at a positive distance.

\bd \cite{ly} [\textbf{Li-Yorke Chaos}]  The system $(X,f)$ is called \textbf{Li-Yorke chaotic} if it has an uncountable scrambled set. \ed

Later it was proved in  \cite{hy} that Devaney chaos implies
Li-Yorke chaos. Actually  Li-Yorke chaos follows from transitivity and the existence of just one periodic
orbit.

\end{section}

\begin{section}{Different Versions of Sensitivity}

\textbf{For any  nonempty, open $U \subset X$ for the system $(X,f)$ we define the seperating times:
$$ N_\delta(U) \quad =  \ \{ n \in \N : diam(f^n(U)) > \delta  \}.$$}

\textbf{$(X, f)$ is sensitive if and only if there exists $\delta > 0$ such
that $N_\delta(U) \neq \emptyset$, for each nonempty, open  $U \subset X$.
}

\bigskip

Recall that $ S \subseteq \N$ is called \textbf{thick} if for each $k \in \N$ there exists $n_k \in \N$ such
that $\{n_k , n_k + 1, . . . , n_k + k\} \subset S$, is called \textbf{cofinite} if it the complement of a finite set in $\N$,   is called \textbf{syndetic} if there exists $M \in \N$ such that $S \cap \{n, n + 1, . . . , n + M\} \neq \emptyset$ for each $n \in \N$, \textbf{piecewise syndetic} if it is the intersection of a syndetic set with a thick set,  \textbf{thickly syndetic} if it has nonempty intersection with every piecewise syndetic
set and is an \textbf{IP-set} if it contains the finite sum of its elements.

\bigskip

Stronger forms of sensitivity have been mainly studied by Christophe Abraham, Gerard Biau and Benoyt Cadre \cite{abc}, Ethan Akin and Sergiy Kolyada \cite{ak}, T. K. Subrahmonian Moothathu \cite{tksm},  Heng Liu, Li Liao, and Lidong Wang \cite{llw}, Xiangdong Ye and Tao Yu \cite{yy}, Wen Huang, Sergiy Kolyada, and Guohua Zhang \cite{hkz}.

\bd The compact system $(X,f)$ is called

\be

\item  \textbf{strongly sensitive} if there exists $\delta > 0$ such that $N_\delta(U)$ is cofinite for any nonempty, open $U \subset X$ \cite{abc}.

\item \textbf{Li-Yorke sensitive} if there exists $\delta > 0$ such that every $x \in X$ is a limit of points $y \in X$ such that the pair $(x, y)$ is proximal but not $\delta-$asymptotic \cite{ak}.

\item  \textbf{thickly sensitive} if there exists $\delta > 0$ such that $N_\delta(U)$ is thick for any nonempty, open $U \subset X$ \cite{tksm}.

\item  \textbf{syndetically sensitive} if there exists $\delta > 0$ such that $N_\delta(U)$ is syndetic for some $m \in \N$ and any nonempty, open $U \subset X$ \cite{tksm}.

\item  \textbf{piecewise syndetically sensitive} if there exists $\delta > 0$ such that $N_\delta(U)$ is piecewise syndetic for  any nonempty, open $U \subset X$ \cite{hkkz, hkz, yy}.

\item  \textbf{thickly syndetically sensitive} if there exists $\delta > 0$ such that $N_\delta(U)$ is thickly syndetic for  any nonempty, open $U \subset X$ \cite{hkkz, hkz, llw, yy}.

\item \textbf{multi-sensitive} if there exists $\delta > 0$ such that $N_\delta(U_i) \neq \emptyset$ for any finite
collection $U_1, \ldots , U_k$ of nonempty, open subsets of $X$ \cite{tksm}.

\item  \textbf{IP sensitive} if there exists $\delta > 0$ such that $N_\delta(U)$ is IP for any nonempty, open $U \subset X$ \cite{hkz, yy}.

\ee \ed

\bex Recall Example \ref{tm}. Here for any open $U \subset [0,1]$, there exists an $n \in \N$ such that $f^n(U) = [0,1]$. Hence, it is easy to see that this system is \emph{strongly sensitive}. \eex

\bex Recall Example \ref{sm}. Here for any pair of open $U, V \subset X$, there exists an $N \in \N$ such that $\sigma^n (U) \cap V \neq \emptyset$ for all $n \geq N$. Thus it is easy to see that this system is \emph{thick sensitive} and \emph{multi-sensitive}. \eex

\bex Recall Example \ref{ms}. This system is sensitive and minimal. It is easy to see that it is \emph{syndetically sensitive}. \eex

From the definition, it follows that

\bigskip

  cofinite sensitivity $\Rightarrow$   thickly syndetical sensitivity $\Rightarrow$ \begin{tabular}{c}

                                                                                      thick sensitivity  \\
                                                                                      syndetical sensitivity \\

                                                                                    \end{tabular}
 $\Rightarrow$   sensitivity.

 \centerline{and}

\centerline{ \small thickly syndetical sensitivity $\Rightarrow$ multi-sensitivity $\Rightarrow$                                                                      thick sensitivity
 $\Rightarrow$   sensitivity $\Leftarrow$  Li-Yorke sensitivity.}

 \bigskip

 We refer to \cite{ak, hkz, llw, tksm, yy} for more examples on different kinds of sensitivity.

 \bigskip

Infact, more can be said

\bt \cite{ak} Let $(X,f)$ be  Li–Yorke sensitive then it is
sensitive. If $(X, f)$ is sensitive and has a fixed point which is its
unique minimal subset then it is Li–Yorke sensitive. \et

\bt \cite{tksm} For interval systems and subshifts of finite type the notions of sensitivity and strong sensitivity are equivalent. \et

\bt \cite{hkz, tksm} A multi-sensitive system is also thickly sensitive. Moreover
if the system is transitive, then the converse also holds. \et

\bt \cite{yy} For minimal systems  the  notions of thick sensitivity, thickly syndetic sensitivity, multi-sensitivity and  IP sensitivity are all equivalent.

        This is also equivalent to the fact that the system is not an almost one-to-one extension of its maximal equicontinuous factor.\et

\br Example \ref{sm} is thickly sensitive, thickly syndetic sensitive, multi-sensitive and IP sensitive. \er

\bigskip

This leads to considering sensitivity for different families of subsets of $\Z_+$.

Let $\mathfrak{P} = \mathfrak{P}(\Z_+)$ be the collection of all subsets of $\Z_+$. $\mathfrak{F} \subset \mathfrak{P}$ is called a \textbf{family} if it is hereditary upwards, i.e. if $F_1 \subset F_2$ and $F_1 \in \mathfrak{F}$ implies $F_2 \in \mathfrak{F}$. A family $\mathfrak{F}$ is \textbf{proper} if it is a proper subset of $\mathfrak{P}$, i.e. it is neither empty nor whole of $\mathfrak{P}$. If a proper family $\mathfrak{F}$ is closed under finite intersection, then $\mathfrak{F}$ is called a filter. We refer to \cite{F} for more details on this.

\bd Let $\mathfrak{F}$ be a family. $(X, f)$ is \textbf{$\mathfrak{F}-$sensitive} if
there is $\delta > 0$ such that for any nonempty, open $U \subset X$, $N_\delta(U) \in \mathfrak{F}$. \ed

\bp \cite{yy} Let $\pi : (X, f) \to (Y, g)$ be an almost-open factor map between two compact systems,
and let $\mathfrak{F}$ be a family. If $(Y, g)$ is $\mathfrak{F}-$sensitive then so is $(X, f)$. \ep

\textbf{$\mathfrak{F}-$sensitivity} for various families $\mathfrak{F}$ have been studied in \cite{yy}.

\end{section}

\begin{section}{Auslander-Yorke Dichotomy}

Recall the \emph{Auslander-Yorke dichotomy} and note the following:

\bp \cite{ay} If the compact system $(X, f)$ is equicontinuous and topologically
transitive, then $(X, f)$ is  minimal, and $T$ is a homeomorphism. \ep

Thus topologically transitive systems if  equicontinuous are minimal and the dichotomy applies. However, not all topologically transitive systems are equicontinuous. These systems are then non minimal. But, a very natural question can be if we can derive Auslander-Yorke dichotomy kind of results for systems not minimal or a stronger version of such a dichotomy.

\medskip

Recall that if the compact  system $(X, f)$ is transitive but not minimal
then the set of non-transitive points is dense.

\bd A topologically transitive system is called \textbf{almost equicontinuous} if there is
at least one equicontinuity point. \ed

We note this dichotomy due to Akin, Auslander and Berg,

\bp \cite{ak} A transitive system is either sensitive or almost equicontinuous. \ep

\medskip

\bd  A system $(X, f)$ is called \textbf{uniformly rigid} if there exists a sequence
$\{n_k\}$ in $\N$ such that the sequence $f^{n_k}$ tends uniformly to the identity map on $X$. \ed

We note this dichotomy due to Glasner and Weiss,

\bp \cite{gw} A topologically transitive system without isolated points which is not
sensitive is uniformly rigid. \ep

\medskip

\bd \cite{hkkz, hkz} The point $x \in X$ is called \textbf{syndetically equicontinuous} for the system $(X,f)$, if there exists a neighborhood $U \ni x$ such that $\{n \in \N: \{f^n\}$ is equicontinuous at $x\}$ is a syndetic set. \ed

Note that the set of syndetically equicontinuous points for $(X,f)$ is a superset of the set of equicontinuity points of $(X,f)$ and that if $(X, f)$  is thickly sensitive then it cannot contain any syndetically equicontinuous points.

\bd \cite{hkkz, hkz} The system $(X, f)$ is \textbf{syndetically
equicontinuous} if every $x \in X$ is syndetically equicontinuous. \ed

\bex Recall Example \ref{ir}. This system is equicontinuous and minimal. It is easy to see that it is syndetically equicontinuous. \eex

\bp \cite{hkkz, hkz} If $(X,f)$ is an almost one-to-one extension of a minimal
equicontinuous system, then it is syndetically equicontinuous.\ep

\bp \cite{hkz} Let $\pi: (X, f) \to (Y, g)$ be an almost one-to-one factor map between
minimal systems. Then $(X,f)$ is syndetically equicontinuous if and only if $(Y, g)$ is also.\ep

We note this dichotomy due to Huang, Kolyada and Zhang,

\bt \cite{hkz} Let $(X, f)$  be a transitive system. Then either $(X, f)$ is thickly sensitive,
 or every transitive point is syndetically equicontinuous. In particular, if $(X, f)$ is minimal,
then it is either thickly sensitive or syndetically equicontinuous. \et

\br Note that in the absence of transitivity the above need not hold. Refer to \cite{hkz} for a counterexample. \er

\bigskip

Some more Auslander-Yorke dichotomy type theorems are obtained in \cite{yy}  involving $\mathfrak{F}-$sensitivity for various families $\mathfrak{F}$.

\bigskip

\textbf{The world awaits more such dichotomies!}
\end{section}


\begin{thebibliography}{99}

\footnotesize

\bibitem{abc}
Christophe Abraham, Gerard Biau and Benoyt Cadre, Chaotic Properties of Mappings on a Probability Space, Journal of Mathematical Analysis and Applications, 266 (2002), 420-431.

\bibitem{aab}
 Ethan Akin, Joseph Auslander, Kenneth Berg, When is a transitive map chaotic?, in: Convergence in ergodic theory and probability, de Gruyter, Berlin-New York, 1996, 25-40.

\bibitem{ak}
Ethan Akin and Sergiy Kolyada, Li–Yorke sensitivity, Nonlinearity 16(4) (2003), 1421-1433.

\bibitem{aus}
Joseph Auslander, Minimal flows and their extensions,  North-Holland Mathematics studies,153 (1988).

\bibitem{ay}
 Joseph Auslander and James  Yorke, Interval maps, factors of maps, and chaos, { Tohoku Mathematical Journal, Second Series,} 32.2, (1980) 177-188.

\bibitem{bb}
 Banks, J. Brooks, J. Cairns, G. Davis, G. and  Stacey, P. On
 Devaney's definition of chaos,  Amer. Math. Monthly, 99,  332-334(1992).

\bibitem{de}
 Robert Devaney,  An introduction to chaotic dynamical systems,
 2nd ed., Addison Wesley, (1989).

\bibitem{F}
{  Hillel Furstenberg,} Disjointness in ergodic theory, minimal sets and a problem in diophantine
approximation, { Mathematical Systems theory,} 1 (1967), 1-49.

\bibitem{gw}
{Eli Glasner and Benjamin Weiss, } Sensitive dependence on initial
conditions, { Nonlinearity,} 6, 1067-1075(1993).


\bibitem{gu}
{ John Guckenheimer, } Sensitive dependence to initial conditions
for one dimensional maps, { Comm. Math. Phys.} 70,
133-160(1979).

\bibitem{hkkz}
Wen Huang, Danylo Khilko, Sergiy Kolyada and Guohua Zhang, Dynamical compactness and sensitivity. J. Differential
Equations 260(9) (2016), 6800-6827.


\bibitem{hkz}
{ Wen Huang, Sergiy Kolyada and Guohua Zhang,} Analogues of Auslander–Yorke theorems for
multi-sensitivity, Ergod. Th. \& Dynam. Sys. (2018), 38, 651-665

\bibitem{hy}
Wen Huang, Xiangdong Ye, Devaney’s chaos or 2-scattering implies Li-Yorke’s chaos, Topology Appl. 117 (2002), No. 3, 259–272.

\bibitem{ly}
Tien-Yien Li and James Yorke,  Period three implies chaos, Am. Math. Month. 82(1975), 985-992.

\bibitem{llw}
Heng Liu, Li Liao, and Lidong Wang, Thickly syndetical sensitivity of topological dynamical
system, Discrete Dyn. Nat. Soc. 2014, Art. ID 583431, 4 pp.

\bibitem{tksm}
T. K. Subrahmonian Moothathu. Stronger forms of sensitivity for dynamical systems. Nonlinearity 20(9)
(2007), 2115-2126.

\bibitem{nk}
{ Anima Nagar  and V. Kannan} Topological Transitivity for Discrete
Dynamical Systems, { Applicable Mathematics In The Golden Age,}
 Narosa Publications(2003), 534-584.

 \bibitem{sil}
{ Stephen Silverman, } On maps with dense orbits and the definition
of chaos, { Rocky Mountain Jour. Math.} 22, 353 - 375(1992).

\bibitem{yy}
X. Ye and Tao Yu. Sensitivity, proximal extension and higher order almost automorphy.  Trans. Amer. Math. Soc. 370 (2018), 3639-3662.


\end{thebibliography}
\end{document}